\documentstyle{article}

\makeatletter
  
  \@addtoreset{equation}{section}
\makeatother

\newtheorem{theorem}{Theorem}
\newtheorem{proposition}{Proposition}
\newtheorem{lemma}{Lemma}
\newtheorem{corollary}{Corollary}
\newtheorem{remark}{Remark}
\newtheorem{example}{Example}

\newfont{\germ}{eufm10}
\newfont{\germlarge}{eufm10}
\newfont{\slsmall}{cmsl8}

\def\clch{\mbox{\sl cl\,ch}\,}
\def\ch{\mbox{\sl ch}\,}
\def\ck#1{{\check #1}}
\def\eps{\epsilon}
\def\et#1{\tilde{e}_{#1}}
\def\ft#1{\tilde{f}_{#1}}
\def\geh{\goth{g}}
\def\goth#1{\mbox{\germ #1}}
\def\ie{{\em i.e., }}
\def\La{\Lambda}
\def\la{\lambda}
\def\li#1{\left[{#1}\right]_+}
\def\ot{\otimes}
\def\ol{\overline}
\def\P{{\cal P}}
\def\Pcl{P_{cl}}
\def\Proof{\noindent{\sl Proof.}\quad}
\def\Q{{\bf Q}}
\def\qed{~\rule{1mm}{2.5mm}}
\def\sln{\widehat{\goth{sl}}_{\,n}}
\def\slt{\widehat{\goth{sl}}_{\,2}}
\def\Uq{U_q(\geh)}
\def\Uqp{U'_q(\geh)}
\def\veps{\varepsilon}
\def\vphi{\varphi}
\def\wt{\mbox{\sl wt}\,}
\def\wts{\mbox{\slsmall wt}\,}
\def\Z{{\bf Z}}
\def\Zn{\Z_{\ge0}}

\begin{document}

\title{ Demazure Modules and Perfect Crystals }

\author{Atsuo Kuniba$^1$, Kailash C. Misra$^2$,\\
        Masato Okado$^3$ and Jun Uchiyama$^4$} 

\date{ \it $^1$Institute of Physics,\\
       \it University of Tokyo,\\
       \it Komaba, Tokyo 153, Japan\\
       \vskip 4mm
       \it $^2$Department of Mathematics,\\
       \it North Carolina State University,\\
       \it Raleigh, NC 27695-8205, USA\\ 
       \vskip 4mm
       \it $^3$Department of Mathematics,\\
       \it The University of Melbourne, \\
       \it Parkville, Victoria, 3052, Australia \thanks{Permanent address: 
       Department of Mathematical Science,
       Faculty of Engineering Science,
       Osaka University, Toyonaka, Osaka 560, Japan.}\\
       \vskip 4mm
       \it $^4$Department of Physics,\\
       \it Rikkyo University,\\
       \it Nishi-Ikebukuro, Tokyo 171, Japan}

\maketitle

\begin{abstract}
\noindent
We give a criterion for the Demazure crystal $B_w(\la)$ defined by
Kashiwara to have a tensor product structure. We study the $\sln$
symmetric tensor case, and see some Demazure characters are expressed
using Kostka-Foulkes polynomials.
\end{abstract} 

\setcounter{section}{-1}

\section{Introduction}

In \cite{S} Sanderson calculated the character of the Demazure module
$V_w(\la)$ of $\slt$ specialized  at $e^{-\delta}=1$ 
($\delta$: the null root) using Littelmann's path \cite{L}. 
To state more precisely, we set
\[
w^{(k)}=\underbrace{r_{1-i}\cdots r_1r_0}_k,
\]
where $r_0,r_1$ are the simple reflections, and $i=0$ ($k$:even), 
$=1$ ($k$:odd). Let $V_{w^{(k)}}(\la)$ be the 
$U(\goth{n}^+)$-module generated by an extremal vector of weight $w^{(k)}\la$
of the irreducible highest weight $U(\slt)$-module of highest weight $\la$. 
Let $V^l$ be the $(l+1)$-dimensional irreducible module of $\goth{sl}_{\,2}$. 
Taking $\la=l\La_0$ for simplicity, her result reads as
\[
\left.\ch V_{w^{(k)}}(l\La_0)\right|_{e^{-\delta}=1}=e^{l\La_i}(\ch V^l)^k.
\]
This suggests that there exists a
tensor product structure in the Demazure module at a combinatorial level.

Now let $\geh$ be any affine Lie algebra. As for the irreducible highest
weight $\Uq$-module $V(\la)$ at $q=0$, we have a combinatorial object 
{\em path}. (This path is different from Littelmann's path.)
It has emerged in the study of solvable lattice models (cf.
\cite{DJKMO1},\cite{DJKMO2}). Its study is accomplished with the aid
of crystal base theory (cf. \cite{KMN1},\cite{KMN2}). Let $B$ be a {\em perfect}
crystal of level $l$. Then, for any dominant integral weight $\la$ of level $l$,
the crystal $B(\la)$ of $V(\la)$ is represented as a set of paths. Roughly
speaking, a path is an element of the semi-infinite tensor product
$\cdots\ot B\ot B$ with some stability condition. On the other hand,
Kashiwara \cite{Ka} gave a recursive formula to obtain the Demazure crystal
$B_w(\la)$ for any symmetrizable Kac-Moody algebra. Denoting the Bruhat
order by $\succ$, it reads as follows.
\begin{eqnarray*}
&&\mbox{If $r_iw\succ w$, then}\\
&&B_{r_iw}(\la)=\bigcup_{n\ge0}\ft{i}^nB_w(\la)\setminus\{0\}.
\end{eqnarray*}
Thus, it seems natural to ask how $B_w(\la)$ is characterized in the set
of paths in the case of affine Lie algebras.

In this article, we present a criterion for the Demazure crystal $B_w(\la)$
to have a tensor product structure. To explain more precisely, let
$w^{(k)}$ ($k=1,2,\cdots$) be an increasing sequence of affine Weyl group
elements with respect to the Bruhat order, let $u_\la$ be a highest
weight vector of $B(\la)$. Then, under the assumptions (I-IV) in section 2.3,
we can show the following isomorphism of crystals.
\[
B_{w^{(k)}}(\la)\simeq u_{\la_j}\ot B_a^{(j,\cdots,j-\kappa+1)}
\ot B^{\ot(j-\kappa)}\qquad\mbox{if}\quad j\ge\kappa.
\]
Here $j,a,\kappa\ge1,B_a^{(j,\cdots,j-\kappa+1)}\in B^{\ot\kappa}$ are determined
when checking the assumptions, and $\la_j=\sigma^j(\la)$ with a Dynkin diagram
automorphism $\sigma$ is defined from the perfect crystal $B$. In Sanderson's case
above, we have $j=k,a=1,\kappa=1,B_a^{(j)}=B$ and $B$ is the crystal of the
$(l+1)$-dimensional irreducible $U_q(\goth{sl}_{\,2})$-module. As a corollary,
we can see
\[
\left.\ch B_{w^{(k)}}(\la)\right|_{e^{-\delta}=1}
=e^{\la_j}(\ch B_a^{(j,\cdots,j-\kappa+1)})(\ch B)^{j-\kappa}
\qquad\mbox{if}\quad j\ge\kappa.
\]

The plan of this article is as follows. In section 1, we summarize the theory
of perfect crystal. We introduce Kashiwara's Demazure crystal $B_w(\la)$,
and state the main result in section 2. We apply our result to the $\sln$
symmetric tensor case in section 3. A relation between Demazure characters
and Kostka-Foulkes polynomials and some other discussions are included in 
section 4.

\bigskip
\noindent{\em Acknowledgement.}\quad
K.C.M. and M.O. would like to thank Omar Foda for discussions, and collaboration
in our earlier work \cite{FMO}. We would like to thank Anatol N. Kirillov for 
suggesting us that in the $\sln$ symmetric tensor case, Demazure characters can be
expressed using Kostka-Foulkes polynomials. We would also like to thank Noriaki
Kawanaka for introducing the book \cite{MP}. M.O. would like to thank Peter
Littelmann for stimulating discussions. K.C.M. is supported in part by 
NSA/MSP Grant No. 96-1-0013. A.K. and M.O. are supported in part by 
Grant-in-Aid for Scientific Research on Priority Areas, the Ministry of 
Education, Science and Culture, Japan. M.O. is supported in part by the Australian
Research Council.

\section{Perfect crystal} 

\subsection{Notation}
We follow the notations of the quantized universal enveloping algebra 
and the crystal base in \cite{KMN1}. In particular, $\{\alpha_i\}_{i\in I}$
is the set of simple roots, $\{h_i\}_{i\in I}$ is the simple coroots, $P$
is the weight lattice and $P_+=\{\la\in P\mid \langle\la,h_i\rangle\ge0
\mbox{ for any }i\}$. $\Uq$ is the quantized universal enveloping algebra of 
an affine Lie algebra $\geh$. $V(\la)$ is the irreducible highest weight module
of highest weight $\la\in P_+$, $u_\la$ is its highest weight vector,
$(L(\la),B(\la))$ is its crystal base.

For the notation of a finite-dimensional representation of $\Uqp$, we follow 
section 3 in \cite{KMN1}. For instance, $\Pcl$ is the classical weight
lattice, $\Uqp$ is the subalgebra of $\Uq$ generated by $e_i,f_i,q^h$
($h\in(\Pcl)^*$) and $\mbox{Mod}^f(\geh,\Pcl)$ is the category of 
finite-dimensional $\Uqp$-modules which have the weight decompositions.
We set $\Pcl^+=\{\la\in\Pcl\mid\langle\la,h_i\rangle\ge0\mbox{ for any }i\}
\simeq\sum\Zn\La_i$ and $(\Pcl^+)_l=\{\la\in\Pcl^+\mid\langle\la,c\rangle=l\}$, 
where $c$ is the canonical central element. Assume $V$ in $\mbox{Mod}^f
(\geh,\Pcl)$ has a crystal base $(L,B)$. For an element $b$ of $B$, we
set $\veps_i(b)=\max\{n\ge0\mid\et{i}^n b\in B\}$, $\veps(b)=\sum_i
\veps_i(b)\La_i$ and $\vphi_i(b)=\max\{n\ge0\mid\ft{i}^n b\in B\}$,
$\vphi(b)=\sum_i\vphi_i(b)\La_i$.

\subsection{Perfect crystal}
In \cite{KMN1} Kang et al. introduced the notion of perfect crystal. 
For the definition of the perfect crystal, see Definition 4.6.1 in \cite{KMN1}.

Let $B$ be a perfect crystal of level $l$. For $\la\in(\Pcl^+)_l$, let 
$b(\la)\in B$ be the element such that $\vphi(b(\la))=\la$. From the definition
of perfect crystal, such a $b(\la)$ exists and is unique. Let $\sigma$ be the
automorphism of $(\Pcl^+)_l$ given by $\sigma\la=\veps(b(\la))$. We set 
$\ol{b}_k=b(\sigma^{k-1}\la)$ and $\la_k=\sigma^k\la$. Then perfectness
assures that we have the isomorphism of crystals
\[
B(\la_{k-1})\simeq B(\la_k)\ot B.
\]
Iterating this isomorphism, we have 
\[
\psi_k:B(\la)\simeq B(\la_k)\ot B^{\ot k}.
\]
Defining the set of paths $\P(\la,B)$ by
\[
\P(\la,B)=\{p=\cdots\ot p(2)\ot p(1)\mid p(j)\in B,p(k)=\ol{b}_k\mbox{ for }k\gg0\},
\]
we see that $B(\la)$ is isomorphic to $\P(\la,B)$. Under this isomorphism, the
highest weight vector in $B(\la)$ corresponds to the path 
$\ol{p}=\cdots\ot \ol{b}_k\ot\cdots\ot \ol{b}_2\ot \ol{b}_1$.
We call $\ol{p}$ the {\em ground-state} path.

\subsection{Actions of $\et{i}$ and $\ft{i}$}
We need to know the actions of $\et{i}$ and $\ft{i}$ on the set of paths 
$\P(\la,B)$. To see this, we consider the following isomorphism. 
\begin{equation}\label{eq:iso1}
\P(\la,B)\simeq B(\la_k)\ot B^{\ot k}.
\end{equation}
For each $p\in\P(\la,B)$, if we take $k$ sufficiently large, we can assume
that $p$ corresponds to $u_{\la_k}\ot p(k)\ot\cdots\ot p(1)$ with $u_{\la_k}$
the highest weight vector of $B(\la_k)$ and $p(n)\in B$ ($n=1,\cdots,k$). 
Then we apply Proposition 2.1.1 in \cite{KN} to see on which component $\et{i}$
or $\ft{i}$ acts. Let us suppose that $\et{i}$ or $\ft{i}$ acts on the $j$-th
component from the right end. If $j<k+1$, we have
\begin{equation}\label{eq:action_e}
\et{i}p=\cdots\ot\et{i}p(j)\ot\cdots\ot p(2)\ot p(1)
\end{equation}
or
\begin{equation}\label{eq:action_f}
\ft{i}p=\cdots\ot\ft{i}p(j)\ot\cdots\ot p(2)\ot p(1).
\end{equation}
If $j=k+1$, we see we should have taken $k$ larger. This happens only for 
$\ft{i}$.

This determination of the component can be rephrased using the notion of
{\em signature}. Let $p\in\P(\la,B)$ correspond to $u_{\la_k}\ot p(k)\ot
\cdots\ot p(1)$ under (\ref{eq:iso1}) as above. With $p(j)$ ($1\le j\le k$), 
we associate 
\begin{eqnarray*}
\eps^{(j)}&=&(\eps^{(j)}_1,\cdots,\eps^{(j)}_m),\\
m&=&\veps_i(p(j))+\vphi_i(p(j)),\\
\eps^{(j)}_a&=&-\quad\bigl(1\le a\le\veps_i(p(j))\bigr),
\quad +\quad\bigl(\veps_i(p(j))<a\le m\bigr).
\end{eqnarray*}
For the highest weight vector $u_{\la_k}$, we take 
\[
\eps^{(k+1)}=(\underbrace{+,\cdots,+}_{\langle\la_k,h_i\rangle}).
\]
We then append these $\eps^{(j)}$'s so that
we have
\[
\eps=(\eps^{(k+1)},\eps^{(k)},\cdots,\eps^{(1)}).
\]
We call it signature of $p$ truncated at the $k$-th position.

Next we consider a sequence of signatures.
\[
\eps=\eta_0,\eta_1,\cdots,\eta_{\max}.
\]
Here $\eta_{j+1}$ is obtained from $\eta_j$ by deleting the leftmost
adjacent $(+,-)$ pair of $\eta_j$. Eventually, we arrive at the 
following signature
\[
\eta_{\max}=(\underbrace{-,\cdots,-}_{n_-},\underbrace{+,\cdots,+}_{n_+}),
\]
with $n_\pm\ge0$. We call it the {\em reduced signature} and denote by
$\ol{\eps}$. 

The component on which $\et{i}$ or $\ft{i}$ acts in (\ref{eq:action_e})
or (\ref{eq:action_f})
reads as follows. If $n_-=0$ (resp. $n_+=0$), we set $\et{i}p=0$ 
(resp. $\ft{i}p=0$). Otherwise, take the rightmost $-$ (resp. leftmost $+$), 
and find the component $\eps^{(j)}$ to which it belonged. Then, this $j$
is the position in (\ref{eq:action_e}) (resp. (\ref{eq:action_f}))
we looked for. Note that if $k$ is large enough, the
position $j$ does not depend on the choice of $k$. 

\begin{remark}
Of course, this signature rule can be applied to the tensor product
of two crystals $B_1\ot B_2$.
\end{remark}

\begin{example}
Let $\geh=\slt$, $B$ be the classical crystal of
the 3-dimensional irreducible representation. Its crystal graph is 
described as follows.
\[
B\hspace{5mm}00\mathop{\rightleftharpoons}_0^1
01\mathop{\rightleftharpoons}_0^1 11
\]
It is known that $B$ is perfect of level 2. We have isomorphisms
$B(\la)\simeq\P(\la,B)$ for $\la=2\La_0,\La_0+\La_1,2\La_1$.
Let $\la=2\La_0$. We see $\la_k=2\La_0$ $(k:\mbox{even})$, $2\La_1$
$(k:\mbox{odd})$. The ground-state path is given by
\[
\ol{p}=\cdots\ot11\ot00\ot11\ot00\ot11.
\]
Consider a path
\[
p=\cdots\ot11\ot01\ot01\ot01\ot00.
\]
The dotted part of $p$ is the same as that of $\ol{p}$. Then the
signature of $p$ truncated at the 5-th position with respect to 
$i=1$ and its reduced signature read as follows.
\begin{eqnarray*}
\eps&=&(++,--,-+,-+,-+,++),\\
\ol{\eps}&=&(\mathop{-}^4\mathop{+}^2\mathop{+}^1\mathop{+}^1).
\end{eqnarray*}
Here the number above each sign shows the component to which it belonged.
Consequently, we have 
\begin{eqnarray*}
\et{1}p&=&\cdots\ot11\ot00\ot01\ot01\ot00,\\
\ft{1}p&=&\cdots\ot11\ot01\ot01\ot11\ot00.
\end{eqnarray*}
\end{example}

\subsection{Weight of a path}
A map $H:B\ot B\rightarrow\Z$ is called {\em energy function} if 
for any $b,b'\in B$ and $i\in I$ such that $\et{i}(b\ot b')\ne0$,
it satisfies
\[
H(\et{i}(b\ot b'))=\left\{\begin{array}{ll}
H(b\ot b')&\quad\mbox{if }i\ne0\\
H(b\ot b')+1&\quad\mbox{if }i=0\mbox{ and }\vphi_0(b)\ge\veps_0(b')\\
H(b\ot b')-1&\quad\mbox{if }i=0\mbox{ and }\vphi_0(b)<\veps_0(b').
\end{array}\right.
\]
It is shown (Proposition 4.5.4 in \cite{KMN1}) that under the 
isomorphism between $B(\la)$ and $\P(\la,B)$, the weight of a path
$p=\cdots\ot p(2)\ot p(1)$ is given by 
\begin{eqnarray}
\wt p&=&\la+\sum_{k=1}^\infty
\left(af(\wt p(k))-af(\wt \ol{b}_k)\right)\nonumber\\
&&\quad-\left(\sum_{k=1}^\infty 
k(H(p(k+1)\ot p(k))-H(\ol{b}_{k+1}\ot \ol{b}_k))\right)\delta. \label{eq:weight}
\end{eqnarray}
For the definition of $af$, along with $cl$, see section 3 of \cite{KMN1}.

\section{Crystal of the Demazure module}

\subsection{Demazure module}
Let $\{r_i\}_{i\in I}$ be the set of simple reflections, and let
$W$ be the Weyl group. For $w\in W$, $l(w)$ denote the length of 
$w$, and $\prec$ denote the Bruhat order on $W$. Let $U_q^+(\geh)$
be the subalgebra of $\Uq$ generated by $e_i$'s. We consider the
irreducible highest weight $\Uq$-module $V(\la)$ ($\la\in\Pcl^+$).
It is known that for $w\in W$, the extremal weight space 
$V(\la)_{w\la}$ is one dimensional. Let $V_w(\la)$ denote the 
$U_q^+(\geh)$-module generated by $V(\la)_{w\la}$. These modules
$V_w(\la)$ ($w\in W$) are called the Demazure modules. They are 
finite-dimensional subspaces of $V(\la)$.

\subsection{Kashiwara's recursive formula}
Let $(L(\la),B(\la))$ be the crystal base of $V(\la)$. In \cite{Ka}
Kashiwara showed that for each $w\in W$, there exists a subset
$B_w(\la)$ of $B(\la)$ such that
\[
\frac{V_w(\la)\cap L(\la)}{V_w(\la)\cap qL(\la)}
=\bigoplus_{b\in B_w(\la)}\Q b.
\]
Furthermore, $B_w(\la)$ has the following recursive property.
\begin{eqnarray}
&&\mbox{If }r_iw\succ w,\mbox{ then}\nonumber\\
&&B_{r_iw}(\la)=\bigcup_{n\ge0}\ft{i}^n B_w(\la)\setminus\{0\}. 
\label{recursive}
\end{eqnarray}
This beautiful formula is essential in our consideration below.

\subsection{Theorem}
We list the assumptions required for our theorem.
Let $\la$ be an element of $(\Pcl^+)_l$, and let $B$
be a classical crystal. Firstly, we assume
\begin{description}
\item[(I)   ]
$B$ is perfect of level $l$.
\end{description}
Therefore, $B(\la)$ is isomorphic to the set of paths
$\P(\la,B)$. Let $\ol{p}=\cdots\ot\ol{b}_2\ot\ol{b}_1$ be
the ground-state path. 

We fix positive integers $d,\kappa$. For a set of elements
$i_a^{(j)}$ in $I$ ($j\ge1,1\le a\le d$), we define
$B_a^{(j)}$ ($j\ge1,0\le a\le d$) by 
\begin{equation} \label{eq:def_B1}
B^{(j)}_0=\{\ol{b}_j\},\hspace{1cm}
B_a^{(j)}=\bigcup_{n\ge0}\ft{i_a^{(j)}}^n B_{a-1}^{(j)}\setminus\{0\}
\quad(a=1,\cdots,d).
\end{equation}
Using these, we next define $B_a^{(j+1,j)}$ ($j\ge1,0\le a\le d$) by
\begin{eqnarray*}
B^{(j+1,j)}_0&=&B_0^{(j+1)}\ot B_d^{(j)},\\
B_a^{(j+1,j)}&=&\bigcup_{n\ge0}\ft{i_a^{(j+1)}}^n B_{a-1}^{(j+1,j)}\setminus\{0\}
\quad(a=1,\cdots,d).
\end{eqnarray*}
Similar definitions continue until we define
\begin{eqnarray*}
B^{(j+\kappa-1,\cdots,j)}_0&=&B_0^{(j+\kappa-1)}\ot B_d^{(j+\kappa-2,\cdots,j)},\\
B_a^{(j+\kappa-1,\cdots,j)}&=&\bigcup_{n\ge0}\ft{i_a^{(j+\kappa-1)}}^n 
B_{a-1}^{(j+\kappa-1,\cdots,j)}\setminus\{0\}
\quad(a=1,\cdots,d).
\end{eqnarray*}
Let us now assume
\begin{description}
\item[(II)  ]
For any $j\ge1$, we have
$B_d^{(j+\kappa-1,\cdots,j)}=B_d^{(j+\kappa-1,\cdots,j+1)}\ot B$.
\end{description}
If $\kappa=1$, the right hand side should be understood as $B$.
We call $\kappa$ a {\em mixing index}. We normally take the minimal one.
The following proposition is sometimes useful. It is proved in the 
same way as in the theorem using Lemma \ref{lemma}. Thus we omit the proof.

\begin{proposition}
If $B_a^{(j+\kappa-1,\cdots,j)}=B_a^{(j+\kappa-1,\cdots,j+1)}\ot B$, then
$B_{a'}^{(j+\kappa-1,\cdots,j)}=B_{a'}^{(j+\kappa-1,\cdots,j+1)}\ot B$ for
$a\le a'\le d$.
\end{proposition}

We introduce the third condition.
\begin{description}
\item[(III) ]
For any $j\ge1$ and $a=1,\cdots,d$, we have
$\langle\la_j,h_{i^{(j)}_a}\rangle\le\veps_{i^{(j)}_a}(b)$
for $b\in B^{(j)}_{a-1}$.
\end{description}
\begin{remark}
Since the automorphism $\sigma$ is of finite order, it suffices to check
the conditions (II),(III) for finite $j$'s.
\end{remark}
For the last assumption, we define $w^{(k)}$ by
\[
w^{(0)}=1,\hspace{1cm}
w^{(k)}=r_{i^{(j)}_a}w^{(k-1)},
\]
where $j$ and $a$ are fixed by $k=(j-1)d+a,j\ge1,1\le a\le d$.
\begin{description}
\item[(IV)  ]
The sequence of Weyl group elements $w^{(0)},w^{(1)},\cdots$
is increasing with respect to the Bruhat order.
\end{description}
This condition is equivalent to $l(w^{(k)})=k$.
The following proposition is convenient to check (IV).

\begin{proposition} \label{pr:Bruhat}
If $\langle w\mu,h_j\rangle>0$ for some $\mu\in\Pcl^+$, then
$r_jw\succ w$.
\end{proposition}

\Proof
{}From $\langle\mu,w^{-1}h_j\rangle>0$, we see $w^{-1}h_j$ is a positive
coroot. This is equivalent to $w^{-1}\alpha_j\in\Sigma_+$, where
$\Sigma_+$ is the set of positive roots. To see this is equivalent to
$r_jw\succ w$, we refer, for example, to proposition 4(i) in section 5.2 (p.411)
of \cite{MP}.
\qed

Now we define a subset $\P^{(k)}(\la,B)$ of $\P(\la,B)$ as follows.
We set $\P^{(0)}(\la,B)=\{\ol{p}\}$.
For $k>0$, we take $j$ and $a$ such that $k=(j-1)d+a,j\ge1,1\le a\le d$,
and set
\begin{equation} \label{eq:def_P}
\P^{(k)}(\la,B)
=\left\{\begin{array}{ll}
\cdots\ot B^{(j+2)}_0\ot B^{(j+1)}_0\ot B^{(j,\cdots,1)}_a
&(j<\kappa)\\
\cdots\ot B^{(j+2)}_0\ot B^{(j+1)}_0\ot B^{(j,\cdots,j-\kappa+1)}_a
\ot B^{\ot(j-\kappa)}
&(j\ge\kappa).
\end{array}\right.
\end{equation}

\begin{theorem}
Under the assumptions (I-IV), we have 
\[
B_{w^{(k)}}(\la)\simeq\P^{(k)}(\la,B).
\]
\end{theorem}

\subsection{Proof of the theorem}
For the proof of the theorem, we need the following simple but
useful lemma.

\begin{lemma} \label{lemma}
Let $B_1,B_2$ be crystals, and let $b_1\in B_1,
b_2\in B_2$. For any $n\ge0$, there exist
$p,q\ge0$ such that 
\[
\ft{i}^p(b_1\ot\et{i}^q b_2)=\ft{i}^n b_1\ot b_2.
\]
\end{lemma}

\Proof
For a reduced signature $\ol{\eps}$,  let $(\ol{\eps})_+$ (resp. 
$(\ol{\eps})_-$) denote the number of $+$'s (resp. $-$'s) in $\ol{\eps}$.
Let $\ol{\eps}_i$ be the reduced signature of
$b_i$ ($i=1,2$). We set
$(\ol{\eps}_1)_+=\alpha,(\ol{\eps}_2)_-=\beta$. 
We divide the proof into two cases: 
(a) $\alpha-n\ge\beta$, 
(b) $\alpha-n<\beta$.

Consider the case (a). In this case, we simply take $p=n,q=0$.
Consider the case (b). In this case, we take $p=\beta-\alpha+2n,
q=\beta-\alpha+n$.
\qed

We proceed to the proof of the theorem.

\noindent{\sl Proof of the theorem.}\quad
In view of assumption (IV), we prove by induction on the length of 
$w^{(k)}$. If $k=0$, \ie $w^{(k)}=1$, the statement
is true. 

Next assume $k>0$ and take $j$ and $a$ such that $k=(j-1)d+a,j\ge1,
1\le a\le d$. From the recursive formula (\ref{recursive})
and $w^{(k-1)}\prec w^{(k)}$, we have 
\[
B_{w^{(k)}}(\la)=\bigcup_{n\ge0}\ft{i^{(j)}_a}^n
B_{w^{(k-1)}}(\la)\setminus\{0\}.
\]
On the other hand, from induction hypothesis, we have
\begin{eqnarray*}
B_{w^{(k-1)}}(\la)&\simeq&
\P^{(k-1)}(\la,B)\\
&=&\hspace{-2pt}\left\{\begin{array}{ll}
\cdots\ot B^{(j+2)}_0\ot B^{(j+1)}_0\ot B^{(j,\cdots,1)}_{a-1}
&(j<\kappa)\\
\cdots\ot B^{(j+2)}_0\ot B^{(j+1)}_0\ot B^{(j,\cdots,j-\kappa+1)}_{a-1}
\ot B^{\ot(j-\kappa)}
&(j\ge\kappa).
\end{array}\right.
\end{eqnarray*}
Note that by definition, this is valid even for $a=1$.
In view of assumption (III), we see the action of $\ft{i^{(j)}_a}$
does not give any effect on the part 
$\cdots\ot B^{(j+2)}_0\ot B^{(j+1)}_0$. 
Thus we ignore this part in the following consideration. Then, in the case
of $j<\kappa$, the proof turns out to be trivial. Assume $j\ge\kappa$,
set $i=i^{(j)}_a,B_1=B^{\ot\kappa},B_2=B^{\ot(j-\kappa)}$,
and take any elements $b_1$ and $b_2$ from 
$B^{(j,\cdots,j-\kappa+1)}_{a-1}\subset B_1$ and $B_2$. 
{}From Lemma \ref{lemma}, for any $n\ge0$, there exist $p,q\ge0$
such that 
\[
\ft{i}^n b_1\ot b_2=\ft{i}^p(b_1\ot\et{i}^q b_2).
\]
Noting that $\et{i}^q b_2\in B^{\ot(j-\kappa)}$, we can conclude
\[
\cdots\ot B^{(j+1)}_0\ot B^{(j,\cdots,j-\kappa+1)}_a\ot B^{\ot(j-\kappa)}
\subset B_{w^{(k)}}(\la).
\]
The other direction of the inclusion is clear.
\qed

\subsection{Application to characters}
Since we have established the weight preserving bijection between
$B_{w^{(k)}}(\la)$ and $\P^{(k)}(\la,B)$, the following proposition
turns out to be an immediate corollary.

\begin{corollary}
Under the assumptions (I-IV), we have 
\[
\ch B_{w^{(k)}}(\la)
=\sum_{p\in\P^{(k)}(\la,B)}e^{\wts p},
\]
where $\wt p$ is given in (\ref{eq:weight}).
\end{corollary}

In view of the tensor product structure of $\P^{(k)}(\la,B)$ (\ref{eq:def_P}),
it would be more interesting to consider a {\em classical character},
which we define by
\begin{eqnarray*}
\clch B_w(\la)&=&\sum_{\mu\in P}\sharp B_w(\la)_\mu e^{cl(\mu)},\\
B_w(\la)_\mu&=&\{b\in B_w(\la)\mid\wt b=\mu\}.
\end{eqnarray*}
Note that for $B(\la)$, the classical character does not 
make sense.

\begin{corollary}
Assume (I-IV). For $k\ge1$, take $j$ and $a$ such that 
$k=(j-1)d+a$, $j\ge1$, $1\le a\le d$. Then,
\[
\clch B_{w^{(k)}}(\la)=\left\{
\begin{array}{ll}
e^{\la_j}\ch B^{(j,\cdots,1)}_a&(j<\kappa)\\
e^{\la_j}(\ch B^{(j,\cdots,j-\kappa+1)}_a)(ch B)^{j-\kappa}&(j\ge\kappa).
\end{array}\right.
\]
\end{corollary}

\Proof 
Let $p$ be an element in $\P^{(k)}(\la,B)$. From (\ref{eq:weight}),
the classical weight of $p$ is given by 
\begin{equation} \label{eq:clch}
\la+\sum_{i=1}^\infty(\wt p(i)-\wt \ol{b}_i).
\end{equation}
Noting that $p(i)=\ol{b}_i$ ($i>j$), $\wt \ol{b}_i=\la_{i-1}-\la_i$
and $\la_0=\la$, (\ref{eq:clch}) reads as
\[
\la_j+\sum_{i=1}^j \wt p(i),
\]
which immediately implies the statement.
\qed

\section{$\widehat{\mbox{\germlarge sl}}_{\,n}$ symmetric tensor case}

In this section, we apply our theorem to the case of symmetric tensor
representations of $U'_q(\sln)$.

\subsection{Preliminaries}
In what follows, $\equiv$ always means $\equiv\mbox{ mod }n$.
Define $\delta_i^{(n)}$ by $\delta_i^{(n)}=1$ ($i\equiv0$),
$=0$ ($i\not\equiv0$). In the case of $\geh=\sln$, we have
$\langle\alpha_i,h_j\rangle=2\delta_{i-j}^{(n)}-\delta_{i-j-1}^{(n)}
-\delta_{i-j+1}^{(n)}$ ($i,j\in I$), where $I=\{0,1,\cdots,n-1\}$.
For our purpose, it is convenient to use the notations $\alpha_i,
h_i,\La_i,r_i$ for $i\in\Z$ by defining $\alpha_i=\alpha_{i'}$,
etc, if $i\equiv i'$. 

Let $B^l$ be the classical crystal of the level $l$ symmetric tensor
representation of $U'_q(\sln)$. As a set, $B^l$ is described as
$B^l=\{(x_0,\cdots,x_{n-1})\in\Zn^n\mid\sum_{i=0}^{n-1}x_i=l\}$.
The actions of $\et{i},\ft{i}$ are defined as follows.
\begin{eqnarray}
\et{i}(x_0,\cdots,x_{n-1})&=&\left\{
\begin{array}{ll}
(x_0,\cdots,x_{i-1}+1,x_i-1,\cdots,x_{n-1})&(i\ne0)\\
(x_0-1,\cdots,x_{n-1}+1)&(i=0)
\end{array}\right.\label{eq:sln_e}\\
\ft{i}(x_0,\cdots,x_{n-1})&=&\left\{
\begin{array}{ll}
(x_0,\cdots,x_{i-1}-1,x_i+1,\cdots,x_{n-1})&(i\ne0)\\
(x_0+1,\cdots,x_{n-1}-1)&(i=0)
\end{array}\right.\label{eq:sln_f}
\end{eqnarray}
If the right hand side contains a negative component, we should understand
it as 0.

Following section 1.2 of \cite{KMN2}, we describe the perfect crystal
structure of $B^l$. (Note that we deal with the case of $(A_{n-1}^{(1)},
B(l\La_1))$ in their notation.) From (\ref{eq:sln_e}) and (\ref{eq:sln_f}),
it is easy to see $\veps\bigl((x_0,\cdots,x_{n-1})\bigr)=\sum_{i=0}^{n-1}
x_i\La_i$, $\vphi\bigl((x_0,\cdots,x_{n-1})\bigr)=\sum_{i=0}^{n-1}
x_i\La_{i+1}$. Setting $\la=\sum_{i=0}^{n-1}m_i\La_i$, we have
$b(\la)=(m_1,\cdots,m_{n-1},m_0)$ and the automorphism $\sigma$ is given
by $\sigma\la=\sum_{i=0}^{n-1}m_i\La_{i-1}$. Thus, the order of $\sigma$ 
is $n$.

It is often convenient to consider automorphisms on $I$ and $B^l$. 
By abuse of notation, we also denote them by $\sigma$. They are defined by
\begin{eqnarray*}
&&\sigma(i)\equiv i-1,\quad 0\le\sigma(i)\le n-1\quad\mbox{ for }i\in I,\\
&&\sigma\bigl((x_0,\cdots,x_{n-1})\bigr)=(x_1,\cdots,x_{n-1},x_0)
\quad\mbox{ for }(x_0,\cdots,x_{n-1})\in B^l.
\end{eqnarray*}
It is easy to see the following properties.
\[
\et{\sigma(i)}\sigma(b)=\sigma(\et{i}b),
\quad\ft{\sigma(i)}\sigma(b)=\sigma(\ft{i}b)
\quad\mbox{ for }b\in B^l.
\]
Using $\sigma$, we have $\ol{b}_j=\sigma^j\bigl((m_0,\cdots,m_{n-1})\bigr)$.

\subsection{Sequence of Weyl group elements}
We define $\li{x}$ by
\[
\li{x}=\left\{
\begin{array}{ll}
\mbox{the largest integer not exceeding }x &(x\ge0)\\
0&(x<0),
\end{array}\right.
\]
and set 
\begin{equation} \label{eq:sln_seq}
w^{(k)}=\left\{
\begin{array}{ll}
1&(k=0)\\
r_{k-1}\cdots r_1r_0&(k>0).
\end{array}\right.
\end{equation}
We are to prove

\begin{proposition}
The sequence of Weyl group elements $w^{(0)},w^{(1)},\cdots$ is 
increasing with respect to the Bruhat order.
\end{proposition}
For the proof, we prepare a lemma.

\begin{lemma}
$w^{(k)}\La_0=\La_0-\sum_{i=0}^{k-1}\li{\frac{i+d}d}\alpha_i$, with
$d=n-1$.
\end{lemma}

\Proof
We prove by induction on $k$. Assuming the statement for $k$ and setting
$\gamma_i=\li{\frac{i+d}d}$, we have
\begin{eqnarray*}
w^{(k+1)}\La_0
&=&w^{(k)}\La_0-\langle w^{(k)}\La_0,h_k\rangle\alpha_k\\
&=&\La_0-\sum_{i=0}^{k-1}\gamma_i\alpha_i\\
&&-\Bigl(\delta_k^{(n)}+\gamma_{k-1}
+\sum_{j\ge1}(\gamma_{k-1-nj}-2\gamma_{k-nj}+\gamma_{k+1-nj})\Bigr)\alpha_k.
\end{eqnarray*}
Therefore, it is sufficient to show
\[
\delta_k^{(n)}+\gamma_{k-1}+\sum_{j\ge1}
(\gamma_{k-1-nj}-2\gamma_{k-nj}+\gamma_{k+1-nj})=\gamma_{k}.
\]
Note that $\gamma_a-\gamma_{a-1}=\theta(a)$,
where $\theta(a)=1$ if $a/d\in\Zn$, $=0$ otherwise. We have
\begin{eqnarray*}
&&\hspace{-5mm}\gamma_k-\gamma_{k-1}
-\sum_{j\ge1}(\gamma_{k-1-nj}-2\gamma_{k-nj}+\gamma_{k+1-nj})\\
&&=\theta(k)+\sum_{j\ge1}\bigl(\theta(k-nj)-\theta(k+1-nj)\bigr)\\
&&=\sum_{j\ge0}\bigl(\theta(k-nj)-\theta(k-nj-d)\bigr)\\
&&=\delta_k^{(n)}.
\end{eqnarray*}
This completes the proof.
\qed

\noindent{\sl Proof of the proposition.}\quad
We take $\La_0$ for $\mu$, and apply Proposition \ref{pr:Bruhat}.
It suffices to show $\langle w^{(k)}\La_0,h_k\rangle>0$.
The left hand side was already calculated to be $\li{\frac{k+d}d}$
in the proof of the previous lemma.
\qed

\subsection{$\la=l\La_0$ case}
We have seen that $B^l$ is perfect of level $l$. We have shown
that the sequence of Weyl group elements $w^{(k)}$ 
(\ref{eq:sln_seq}) is increasing with respect to the Bruhat
order. Thus, taking $d=n-1$ and $i_a^{(j)}\equiv a-j$
($0\le i_a^{(j)}\le n-1$), the assumptions (I) and (IV) are
already satisfied. Let us take $\la=l\La_0$. Then, we have
$\la_j=l\La_{-j}$ and the ground-state path is given by
$\ol{p}=\cdots\ot\ol{b}_2\ot\ol{b}_1$, 
$\ol{b}_j=(0,\cdots,0,\stackrel{i}l,0,\cdots,0)$ with
$i\equiv-j$ ($0\le i\le n-1$). Let us describe $B_a^{(j)}$. Thanks to
the automorphism $\sigma$, it is sufficient to consider the case 
of $j=n$. From (\ref{eq:def_B1}) and $i_a^{(n)}=a$, we easily 
have
\[
B_a^{(n)}=\{(x_0,\cdots,x_{n-1})\in B^l\mid\sum_{i=0}^a x_i=l\}.
\]
Checking the assumption (III) is also trivial. Therefore, we have
the following proposition.

\begin{proposition}
For the case of $\la=l\La_0$, take $d=n-1$ and $i_a^{(j)}\equiv
a-j$ $(0\le i_a^{(j)}\le n-1)$. Then, the assumptions (II) and 
(III) are satisfied with $\kappa=1$.
\end{proposition}

\subsection{$\la$ arbitrary case}
We show even if $\la$ is any element in $(\Pcl^+)_l$, we have 
$\kappa=2$ with the same choice of $d$ and $i_a^{(j)}$ as in 
the previous subsection. For $\la=\sum_{i=0}^{n-1}m_i\La_i$,
define a subset $B_\la^l$ of $B^l$ by
\[
B_\la^l=\{(x_0,\cdots,x_{n-1})\in B^l\mid
x_0+\cdots+x_{i-1}\le m_0+\cdots+m_{i-1}\mbox{ for }1\le i\le n-1\}.
\]
We now prepare a lemma.

\begin{lemma}
Let $\la=\sum_{i=0}^{n-1}m_i\La_i$, $\sum_{i=0}^{n-1}m_i=l_1$.
For any $b_1\in B_\la^{l_1}$ and $b_2\in B^{l_2}$, there exist an
element $\ck{b}_2\in B^{l_2}$ and integers $p_1,\cdots,p_{n-1}\ge0$
such that $\ft{n-1}^{p_{n-1}}\cdots\ft{1}^{p_1}(\ck{b}_1\ot\ck{b}_2)
=b_1\ot b_2$, $z_i\le x_{i-1}$ $(i=1,\cdots,n-1)$, where we set
$\ck{b}_1=(m_0,\cdots,m_{n-1})$, $\ck{b}_2=(z_0,\cdots,z_{n-1})$,
$b_1=(x_0,\cdots,x_{n-1})$.
\end{lemma}

\Proof
We prove by induction on $n$. If $n=1$, the statement is trivial. 
Now assume the statement is valid when $n-1$, and set 
$b_2=(y_0,\cdots,y_{n-1})$. We divide into two cases:(a) $x_{n-2}
>y_{n-1}$, (b) $x_{n-2}\le y_{n-1}$.

Consider the case (a). Assume $(x_0,\cdots,x_{n-1})\in B_\la^{l_1}$,
then we have $x_{n-1}\ge m_{n-1}$. Consider elements
$(x_0,\cdots,x_{n-3},x_{n-2}+x_{n-1}-m_{n-1},m_{n-1})\in B_\la^{l_1}$,
$(y_0,\cdots,y_{n-1})\in B^{l_2}$. Ignoring the last component and
using the induction hypothesis, we see there exist $z_0,\cdots,z_{n-2}$
and $p_1,\cdots,p_{n-2}\ge0$ such that $(z_0,\cdots,z_{n-2},y_{n-1})
\in B^l$, and
\begin{eqnarray*}
&&\ft{n-2}^{p_{n-2}}\cdots\ft{1}^{p_1}
\bigl((m_0,\cdots,m_{n-1})\ot(z_0,\cdots,z_{n-2},y_{n-1})\bigr)\\
&&\hspace{5mm}=(x_0,\cdots,x_{n-3},x_{n-2}+x_{n-1}-m_{n-1},m_{n-1})
\ot(y_0,\cdots,y_{n-1}).
\end{eqnarray*}
Apply $\ft{n-1}^{x_{n-1}-m_{n-1}}$ further, then we get
$(x_0,\cdots,x_{n-1})\ot(y_0,\cdots,y_{n-1})$. Checking another
condition is easy.

We move to (b). The proof goes similarly. Take 
$\ck{b}_2=(z_0,\cdots,z_{n-2},x_{n-2})$ ($z_0,\cdots,z_{n-2}$
are determined from the induction hypothesis) and $p_{n-1}=
x_{n-1}-m_{n-1}+y_{n-1}-x_{n-2}$.
\qed

\begin{proposition}
For any $\la\in(\Pcl^+)_l$, take $d=n-1$ and $i_a^{(j)}\equiv a-j$
$(0\le i_a^{(j)}\le n-1)$. Then, the assumptions (II) and (III) are
satisfied with $\kappa=2$.
\end{proposition}

\Proof
We again reduce the proof for all $j$ to the $j=n$ case by the
automorphism $\sigma$. Set $\la=\sum_{i=0}^{n-1}m_i\La_i$. Note
that $i_a^{(n)}=a$. Firstly, we claim that $B_a^{(n)}$ consists of the 
elements $(x_0,\cdots,x_{n-1})\in B^l$ satisfying
\begin{eqnarray*}
&x_0+\cdots+x_{i-1}\le m_0+\cdots+m_{i-1}\quad&(1\le i\le a),\\
&x_i=m_i&(a+1\le i\le n-1).
\end{eqnarray*}
This is easily shown by induction on $a$. Thus, we have
$\langle\la_n,h_a\rangle=\veps_a(B_{a-1}^{(n)})=m_a$, and (III) is checked.

It remains to show $B_d^{(n,n-1)}=B_d^{(n)}\ot B^l$. Noting that
$B_d^{(n)}=B_\la^l$, it is equivalent to check
\begin{equation} \label{eq:to_check}
\bigcup_{p_1,\cdots,p_{n-1}\ge0}
\ft{n-1}^{p_{n-1}}\cdots\ft{1}^{p_1}
\bigl(\{\ol{b}_n\}\ot\sigma^{-1}(B_\la^l)\bigr)\setminus\{0\}
=B_\la^l\ot B^l.
\end{equation}
Take any $b_1=(x_0,\cdots,x_{n-1})\in B_\la^l$ and $b_2\in B^l$. 
{}From Lemma, there exist $\ck{b}_2=(z_0,\cdots,z_{n-1})\in B^l$
and $p_1,\cdots,p_{n-1}\ge0$ such that
$\ft{n-1}^{p_{n-1}}\cdots\ft{1}^{p_1}(\ol{b}_n\ot\ck{b}_2)$
$=b_1\ot b_2$, $z_i\le x_{i-1}$ ($i=1,\cdots,n-1$).
Noting that $z_1+\cdots+z_i\le x_0+\cdots+x_{i-1}\le m_0+\cdots+m_{i-1}$,
we see $\ck{b}_2\in\sigma^{-1}(B_\la^l)$. Thus, we have shown the 
inclusion $\supset$ in (\ref{eq:to_check}).

The other direction of the inclusion is clear.
\qed

\section{Discussion}

We explain the suggestion by Kirillov. Consider the case of $\geh=\sln,
B=B^l$. We retain the notations in section 3. If $L$ is divisible by $n$,
we already know
\[
B_{w^{(Ld)}}(l\La_0)\simeq u_{l\La_0}\ot(B^l)^{\ot L}.
\]
Therefore, $B_{w^{(Ld)}}(l\La_0)$ is invariant under the action of 
$\et{i},\ft{i}$ ($i\ne0$). On the other hand, it has been shown recently
that the Kostka-Foulkes polynomial $K_{\la(l^L)}(q)$ has the following
expression (cf. \cite{DF},\cite{D},\cite{NY}).
\[
K_{\la(l^L)}(q)=\sum q^{\sum_{j=1}^{L-1}jH(b_{j+1}\ot b_j)},
\]
where the sum is over all highest weight vectors $b_L\ot\cdots\ot b_1$
with respect to $\et{i}$ ($i\ne0$) with highest weight $\sum_{i=1}^{n-1}
(\la_i-\la_{i+1})\ol{\La}_i$. Now Kirillov's suggestion reads as
\begin{eqnarray*}
&&\left(e^{-l\La_0}\ch B_{w^{(Ld)}}(l\La_0)\right)(z_1,\cdots,z_{n-1};q)\\
&&\hspace{8mm}=\sum_\la q^{-E_0}
\left.K_{\la(l^L)}(q)s_\la(x_1,\cdots,x_n)\right|_{x_1\cdots x_n=1},
\end{eqnarray*}
where $\la$ runs over all partitions of $lL$ of depth less or equal than
$n$, $E_0=\frac{lL}2(\frac{L}n-1)$, $q=e^{-\delta}$, $z_i=e^{-\alpha_i}
=x_{i+1}/x_i$ ($i=1,\cdots,n-1$), and $s_\la$ is the Schur function.
Note that an expression of $K_{\la\mu}(q)$ using Gaussian polynomials
is known \cite{Ki}.

In this article, we have presented a criterion for the Demazure crystal 
to have a tensor product structure. One may ask if there are other cases
than the $\sln$ symmetric tensor case which have a similar property. We claim
that if $\geh$ is non-exceptional and $\la=l\La_0$ (plus some other $\la$), 
we can find a sequence of Weyl group elements which satisfies the assumptions 
(II-III) with $\kappa=1$. We would like to report on that in a near future.

We are not only interested in the classical character, but also in the full
character. One of the advantages to relate the Demazure crystal to a set
of paths is that we can utilize the results on so-called 1D sums in 
solvable lattice models. In \cite{FMO} we have seen that for $\slt$ case, 
all characters of the Demazure modules are expressed in terms of 
$q$-multinomial coefficients. We also would like to deal with this subject for 
any non-exceptional affine Lie algebra (though the level will be 1) in
another publication.

\end{document}